\RequirePackage{ifpdf}
\ifpdf 
\documentclass[pdftex]{sigma}
\else
\documentclass{sigma}
\fi

\DeclareMathOperator{\Hom}{Hom} \DeclareMathOperator{\kernel}{Ker}
\DeclareMathOperator{\sign}{sign}
\DeclareMathOperator{\trace}{trace}
\DeclareMathOperator{\Trace}{Tr}
\DeclareMathOperator{\Wres}{Res} 
\DeclareMathOperator{\wres}{res} 

\begin{document}


\renewcommand{\PaperNumber}{104}

\FirstPageHeading

\renewcommand{\thefootnote}{$\star$}

\ShortArticleName{Conformal Invariants for Manifolds with
Boundary}

\ArticleName{Some Conformal Invariants from the Noncommutative
Residue for Manifolds with Boundary\footnote{This paper is a
contribution to the Proceedings of the 2007 Midwest Geometry
Conference in honor of Thomas~P.\ Branson. The full collection is
available at
\href{http://www.emis.de/journals/SIGMA/MGC2007.html}{http://www.emis.de/journals/SIGMA/MGC2007.html}}}

\Author{William J. UGALDE}

\AuthorNameForHeading{W.J. Ugalde}

\Address{Escuela de Matem\'atica, Universidad de Costa Rica,
C\'odigo postal 2060 San Jos\'e, Costa Rica}
\Email{\href{mailto:william.ugalde@ucr.ac.cr}{william.ugalde@ucr.ac.cr}}
\URLaddress{\url{http://www2.emate.ucr.ac.cr/~ugalde/}}

\ArticleDates{Received August 06, 2007, in f\/inal form October
31, 2007; Published online November 07, 2007}

\Abstract{We review previous work of Alain Connes, and its
extension by the author, on some conformal invariants obtained
from the noncommutative residue on even dimensional compact
manifolds without boundary.  Inspired by recent work of Yong Wang,
we also address possible generalizations of these conformal
invariants to the setting of compact manifolds with boundary.}


\Classification{53A30}

\section{Introduction}

There is one particular aspect of noncommutative geometry that has
historically received less attention than other of its subjects;
the use of its machinery to obtain conformal invariants
(associated to the underlying manifold).   The motivating example
in this venue, is a conformal invariant in dimension~4
(Connes~\cite{ConQuan}) and its extension to higher order even
dimensional mani\-folds by the author~\cite{MioCMP}. The main idea
lies on Theorem~IV.4.2.c of Connes~\cite{Koran}.  This theorem
states that the oriented conformal structure of a compact
even-dimensional smooth manifold is uniquely determined by the
Fredholm module $(\mathcal{H},F,\gamma)$ of Connes, Sullivan and
Teleman~\cite{ConnesST}; via the noncommutative residue $\Wres$ of
Adler, Manin, Guillemin, and Wodzicki~\cite{Adler, Manin,
Guillemin, Wodzicki}.

In Section~2 of \cite{ConQuan} Connes uses his quantized calculus
to f\/ind a conformal invariant in the 4-dimensional case.  A
central part of the explicit computation of this conformal
invariant  is the study of a trilinear functional on smooth
functions over the manifold $M^4$ given by the relation
\[
\tau(f_0,f_1,f_2) = \Wres(f_0[F,f_1][F,f_2]).
\]
Here $F$ is a pseudodif\/ferential operator of order 0 acting on
2-forms over $M^4.$  This conformal invariant computed by Connes
in the $4$-dimensional case is a natural bilinear dif\/ferential
functional  of order $4$  acting on $C^\infty(M^4).$  In
\cite{ConQuan} it is denoted $\Omega$ and in these notes it is
denoted $B_4 \,d x.$

This bilinear functional is symmetric, $B_4(f_1,f_2) =
B_4(f_2,f_1),$ and conformally invariant, in the sense that
$\widehat{B_4}(f_1,f_2) = e^{-4\eta}B_4(f_1,f_2)$ for a conformal
change of the metric $\widehat{g} = e^{2\eta}g.$ It is also
uniquely determined by the relation:
\[
\tau(f_0,f_1,f_2) = \int_M f_0 B_4(f_1,f_2)\,d  x,\qquad \forall\,
f_i \in C^\infty(M^4).
\]
Furthermore, in the 4-dimensional case, Connes has also shown that
the Paneitz operator \cite{Paneitz} (\emph{critical} {\sf GJMS}
for $n=4$ \cite{GJMS}), can be derived from $B_4$ by the relation
\[
\int_M B_4(f_1,f_2)\,d  x = \frac{1}{2} \int_M f_1 P_4(f_2)\,d x.
\]

Aiming to extend the work of Connes to even dimensional manifolds,
in \cite{Ugalde} we have proved the following two results:

\medskip

\noindent \textbf{Theorem~1 of~\cite{Ugalde}.} \textit{Let $M$ be
an $n$-dimensional compact conformal manifold without boundary.
Let $S$ be a pseudodifferential operator of order 0 acting on
sections of a vector bundle over $M$ such that $S^2 f_1 = f_1 S^2
$ and the pseudodifferential operator $P=[S,f_1][S,f_2]$ is
conformally invariant for any $f_i \in C^\infty(M).$ Then there
exists a unique, symmetric, bilinear, differential functional
$B_{n,S}$ of order $n$ conformally invariant in the sense that
$\widehat{B_{n,S}}(f_1,f_2) = e^{-n\eta}B_{n,S}(f_1,f_2)$, for
$\widehat{g} = e^{2\eta}g$, and such that
\[
\Wres(f_0[S,f_1][S,f_2]) = \int_M f_0 B_{n,S}(f_1,f_2)\, d x
\]
for all $f_i \in C^\infty(M).$}

\medskip

A particular case of the above result occurs when one works with
even-dimensional manifolds for then, it makes sense to consider
the Fredholm module $(\mathcal{H},F)$ associated to $M.$  The
operator $F$ has the property $F^2=1$ and in general $[F,f]
\not=0$ for $f \in C^\infty(M).$  Taking $S$ as $F$ in the
previous theorem one has:

\medskip

\noindent \textbf{Theorem~2 of \cite{Ugalde}.} \textit{Let $M$ be
a compact conformal manifold without boundary of even dimension
$n$ and let $(\mathcal{H},F)$ be the Fredholm module associated to
$M$ by A.~Connes~{\rm \cite{ConQuan}}.  Then, by taking $S=F$ in
Theorem~{\rm 1} of {\rm \cite{Ugalde}} there is a unique, symmetric, and
conformally invariant $n$-differential form $B_n=B_{n,F}$ such
that
\[
\Wres(f_0[F,f_1][F,f_2]) = \int_M f_0 B_{n}(f_1,f_2)\,d x
\]
for all $f_i \in C^\infty(M).$}

\medskip

These results are based on the study of the formula for the total
symbol $\sigma(P_1 P_2)$ of the product of two
pseudodif\/ferential operators, in the particular case in which
one of them is a~multiplication operator. The following is the
main result of \cite{MioCMP}:

\begin{theorem}
Let $M$ be a compact conformal manifold without boundary of even
dimension $n$ and let $(\mathcal{H},F)$ be the Fredholm module
associated to $M$ by Connes {\rm \cite{ConQuan}}. Let $P_n$ be the
differential operator given by the relation
\[
\int_M B_n(f,h) \,d x = \int_M f P_n(h) \,d x
\]
for all $f,h \in  C^\infty(M).$  Then,
\begin{itemize}\itemsep=0pt
\item[i)] $P_n$ is formally selfadjoint; \item[ii)] $P_n$ is
conformally invariant in the sense $\widehat{P_n}(h) = e^{-n\eta}
P_n(h),$ if $\widehat g = e^{2 \eta}g$; \item[iii)] $P_n$ is
expressible universally as polynomial in the components of
$\nabla$ (the covariant derivative) and $R$ (the curvature tensor)
with coefficients rational in $n.$ \item[iv)] $P_n(h) = c_n
\Delta^{n/2}(h) +$ ``lower order terms'', with $c_n$ a universal
constant; \item[v)] $P_n$ has the form $\delta S_n d$ where $S_n$
is an operator on 1-forms given as a constant multiple of
$\Delta^{n/2-1} + $ ``lower order terms'' or $(d\delta)^{n/2-1} +$
``lower order terms''; \item[vi)] $P_n$ and $B_n$ are related by:
\[
P_n(f h) - f P_n(h) - h P_n(f) = -2B_n(f,h).
\]
\end{itemize}
\end{theorem}

Because the critical {\sf GJMS} operator and the operator $P_n$
coincide in the f\/lat case and share the same conformal behavior
we have
\begin{proposition}
In the even dimensional case, inside the conformally flat class of
metrics, the critical {\sf GJMS} operator and the operator $P_n$
coincide up to a constant multiple.
\end{proposition}
%

\subsection{Yong Wang's work}
Based on the work of \cite{mio1}, in \cite{Wang} and \cite{Wang2}
Y. Wang proposes to extend to the case of manifolds with boundary,
the work of Connes in Section~2 of \cite{ConQuan}.  He is the
f\/irst to suggest the replacement of  the  usual noncommutative
residue by the noncommutative residue of
Fedosov--Golse--Leichtnam--Schrohe \cite{Fedosovetal}, acting on
Boutet de Monvel's algebra \cite{Boutet}.

In Section 3 of \cite{Wang}, Y.~Wang considers a compact
$n$-dimensional manifold $X$ with boundary~$Y$ and its double
manifold $\widetilde{X}= X\cup_Y X.$  For a vector bundle $E$ over
$\widetilde{X}$ and a pseudodif\/ferential operator $S$ with the
transmission property and of order $0$ acting on sections of~$E$,
the operator~$\widetilde {P}$ is def\/ined as the composition
\[
\widetilde {P} :=
\begin{pmatrix}
\pi^+ f_0 & 0
\\
0&0
\end{pmatrix}
\left[
\begin{pmatrix}
\pi^+ S & 0
\\
0&0
\end{pmatrix},
\begin{pmatrix}
\pi^+ f_1 & 0
\\
0&0
\end{pmatrix}
\right] \left[
\begin{pmatrix}
\pi^+ S & 0
\\
0&0
\end{pmatrix},
\begin{pmatrix}
\pi^+ f_2 & 0
\\
0&0
\end{pmatrix}
\right].
\]
It is then observed that $\widetilde{P} = \pi^+(f_0[S,f_1][S,f_2])
+ G$ for some singular Green operator $G$ with singular Green
symbol $b.$ See \cite{Wang} for the corresponding def\/initions.
Based on this decomposition of $\widetilde{P}$ and the
def\/inition of Fedosov et al.~of~$\overline{\Wres}$, Wang
def\/ines $\Omega_{n,S}$ and $\Omega_{n-1,S}$ via
\[
\Omega_{n,S}(f_1,f_2) = \overline{\wres}
\bigl(([S,\overline{f}_1][S,\overline{f}_2])|_X\bigr) \qquad
\hbox{and} \qquad f_0|_Y\Omega_{n-1,S}(f_1,f_2) = 2\pi \wres_{x'}
\trace(b),
\]
with $\overline{f}_i$ a smooth extension of $f_i$ to
$\widetilde{X}.$  Also $\overline{\wres}$ is the density
corresponding to $\overline{\Wres}$ (similar to $\wres$ for
$\Wres$) and $\wres_{x'}$ the noncommutative residue density for
the manifold $Y.$

It is possible to verify that $\widetilde{P}$ satisf\/ies the
transmission property.  In this way, forgetting about any
conformal invariance property, in \cite{Wang} Wang found a
generalization of the relation
\[
\Wres(f_0[S,f_1][S,f_2]) = \int_M f_0 \Omega_{n,S}(f_1,f_2)
=\int_M f_0 B_{n,S}(f_1,f_2)\, d x, \qquad\forall\,f_i \in
C^\infty(X\cup Y)
\]
in Theorem~1 of \cite{Ugalde}.

The main idea of Y. Wang relies on the use of the double manifold.
Topologically, the double manifold of a given compact oriented
manifold with boundary makes perfect sense. At the level of smooth
manifold with a given Riemannian structure more work is needed to
make sense of a~double manifold.

In Section~4 of \cite{Wang} the dimension is taken as even and the
metric on $X$ has a product structure near the boundary: $g^X =
g^{\partial X} + d x^n.$ On $\widetilde{X}$ the metric
$\widetilde{g}$ is taken as $\widetilde{g}=g$ on both copies of
$X.$ Then $\Omega_n$ is def\/ined as
$\Omega_{n,F}(\overline{f_1},\overline{f_2})|_X$ where
$(\mathcal{H},F)$ is the Fredholm module associated to
$(\widetilde{X},\widetilde{g}),$ and $\overline{f_i}$ is an
extension of $f_i$ to $\widetilde{M}.$

In \cite{Wang2} the even-dimensional Riemannian metric in
consideration has the particular form $g^X = 1/(h(x_n))
g^{\partial X} + d x_n^2$ on a collar neighborhood $U$ of
$\partial X.$ Here $h$ is the restriction to $[0,1)$ of a smooth
function $\widetilde{h}$ on $(-\varepsilon,1)$ for some
$\varepsilon>0$ such that $h(0)=1$ and $h(x_n)>0.$ A metric
$\hat{g}$ is associated to the double manifold $\widetilde{X}$ in
the following way.  Given $U$ and $\widetilde{h}$ as before, there
is a metric $\hat{g}$ on $\widetilde{X}$ with the form
$\hat{g}^{\widetilde{X}} = 1/(\widetilde{h}(x_n))g^{\partial X} +
dx_n^2$ on $U\cup_{\partial X}\partial X\times(-\varepsilon,0]$
and such that $\hat{g}|_X = g.$

Next, with $(\mathcal{H},F_{\hat{g}})$ the Fredholm module
associated to $(\widetilde{X},\hat{g}),$ $\Omega_n$ and
$\Omega_{n-1}$ are def\/ined via the relation
\[
\overline{\Wres} \bigl(\pi^+ f_0 [\pi^+ F_{\hat{g}},f_1][\pi^+
F_{\hat{g}},f_2]\bigr) = \int_X f_0\Omega_n(f_1,f_2)(\hat{g}) +
\int_{\partial M} f_0|_{\partial X} \Omega_{n-1}(\hat{g}).
\]

The described settings used in \cite{Wang} and \cite{Wang2} have
the following limitations:  if we conformally rescale the metric
in $X$ then, the new metric $e^{2\eta}g$ is not anymore of the
specif\/ic requested form near the boundary.  How to def\/ine then
the objects in question in terms of this new metric?  That is to
say, what is the def\/inition of $\Omega_n(e^{2\eta}g)$ and how to
compare it with $\Omega_n(g)$?  Based on the idea of replacing
$\Wres$ with $\overline{\Wres},$ and inspired by  the work of
Wang, we propose the approach in this work to the problem of
extending the results in Section~2 of \cite{ConQuan} to manifolds
with boundary.

\subsection{Contents}

We f\/irst review the construction of the even Fredholm module
$(\mathcal{H},F,\gamma)$ over the commutative algebra $\mathcal{A}
= C^\infty(M)$ (trivially an involutive algebra over $\mathbb{C}$)
of smooth functions over a compact oriented manifold $M$ without
boundary.  We give special attention to its conformal properties.
Then we review the statement of Connes about recovering the
conformal structure from this Fredholm module and a recent
characterization on the subject by B\"ar.  Next, we move to the
setting of manifolds with boundary. Aiming to extend previous
work, we brief\/ly review the noncommutative residue for manifolds
with boundary and Boutet de Monvel calculus according to our
needs.  Last, we present a couple of results that extend to the
even dimensional case Theorems~1~and~2 in \cite{Ugalde} to the
following setting of manifolds with boundary: $M$ is a compact
manifold with boundary $\partial M$ such that $M$ is embedded in a
compact oriented manifold $\widetilde M$ without boundary.
Further we assume Riemannian structures
$(\widetilde{M},\widetilde{g})$ and $(M,g)$ such that $g$
coincides with $\widetilde{g}$ restricted to $M.$ The results are
Theorems~\ref{Theorem1boundary} and~\ref{Theorem2boundary}
respectively.

\subsection{Other possibilities}

For a Riemann surface $M$, a map $f = (f^i)$ from $M$ to
$\mathbb{R}^2$ and metric $g_{i j}(x)$ on $M$, the 2-di\-mensional
Polyakov action \cite{Pol} is given by
\[
I(f) = \frac{1}{2\pi} \int_M g_{i j} \,d f^i \wedge \star d f^j.
\]

By considering instead of $d f$ its quantized version $[F,f],$
Connes \cite{ConQuan} quantized the Polyakov action as a Dixmier
trace:
\begin{align*}
\frac{1}{2 \pi} \int_M g_{i j} d f^i \wedge \star d f^j =
-\frac{1}{2} \Trace_\omega\bigl(g_{i j}[F,f^i][F,f^j]\bigr).
\end{align*}
Connes' trace theorem \cite{ConAction} states that the Dixmier
trace and the noncommutative residue of an elliptic
pseudodif\/ferential operator of order $-n$ on an $n$-dimensional
manifold $M$ are proportional by a factor of $n(2\pi)^n.$  In the
2-dimensional case the factor is $8\pi^2$ and so, the quantized
Polyakov action can be written as
\begin{align*}
-16\pi^2 I=\Wres \bigl(g_{i j} [F,f^i][F,f^j]\bigr).
\end{align*}
This quantized Polyakov action makes sense in the general even
dimensional case.

Although this generalization of the Polyakov action motivates the
particular form of the functional
 $\Wres(f_0 [F,f_1] [F,f_2]),$ the same Fredholm module yields other
functionals in dimensions greater than~4. For instance, for
dimension~6 one could also consider
\[
\Wres(f_0[F,f_1][F,f_2][F,f_3]) = \int_M f_0\, T(f_1,f_2,f_3)
\,d^6x,
\]
which is a Hochschild 3-cocycle and the trilinear expression
$T(f_1,f_2,f_3)$ is conformally invariant. In greater $2l$
dimensions,
\[
\Wres(f_0[F,f_1]\cdots[F,f_{l}]) =  \int_M f_0\,
C(f_1,\dots,f_{l}) \,d^{2l}x,
\]
invites to study the role of the conformal invariant
$C(f_1,\dots,f_{l}).$

\subsection{Further directions}

The possibility of obtaining from the expression
\begin{gather*}
\int_M B_n(f_1,f_2\bigr) \, d x +  2\pi \int_{\partial M} \partial B_n(f_1|_{\partial M},f_2|_{\partial M}) \, d x'\\
\qquad {} = \int_M f_1 P_n f_2\,d x  + \int_{\partial M}
f_1|_{\partial M} P'_{n-1} f_2|_{\partial M} \, d x'
\end{gather*}
conformally covariant dif\/ferential operators $P_n$ and
$P'_{n-1}$ acting on $M$ and $\partial M$ respectively (in a way
similar to the boundaryless case), is the motivating force behind
this project. We hope to report on that in the near future.

One more possibility is to study a Riemannian manifold  with a
particular metric structure near the boundary, in such a way that
it makes sense to consider its double manifold and at the same
time, study conformal variations of the metric.  One can ask what
sort of specif\/ic objects are to be found using the ideas
presented here in such a particular situation.    A seemingly
promising case is that of manifolds with totally geodesic
boundaries, for which the double manifold is natural to be
considered.

\section{The Fredholm module for a conformal manifold}

Following Def\/inition IV.4.1~\cite{Koran}, an even Fredholm
module $(\mathcal{H},F,\gamma)$ is given by
\begin{itemize}\itemsep=0pt
\item An involutive algebra $\mathcal{A}$ (over $\mathbb{C}$)
together with a Hilbert space $\mathcal{H}$ and an involutive
representation $\pi$ of $\mathcal{A}$ in $\mathcal{H}.$ \item An
operator $F$ on $\mathcal{H}$ such that $F=F^*,$ $F^2=1,$ and
$[F,\pi(a)]$ is a compact operator $\forall\,a\in \mathcal{A}.$
\item A $\mathbb{Z}/2$ grading $\gamma,$ $\gamma=\gamma^*,$
$\gamma^2=1$ of $\mathcal{H}$ such that $\gamma \pi(a) = \pi(a)
\gamma, \, \forall\,a\in \mathcal{A},$ and $\gamma F = - F
\gamma.$
\end{itemize}

In the case of a manifold without boundary, the very f\/irst
ingredient in $(\mathcal{H},F,\gamma)$ is the involutive algebra
$\mathcal{A} = C^\infty(M)$ where we allow complex values.   The
fact that, for an oriented  Riemannian manifold $M$ (with or
without boundary) and of even dimension $n,$ the restriction of
the Hodge star operator to middle-dimension forms is conformally
invariant, is central to what follows.  We consider the vector
bundle $\Omega_{\mathbb{C}}^{n/2}(M)$ of complex middle-dimension
forms.  We drop the subscript $\mathbb{C}$ from now on. In this
way, for an $n$-dimensional oriented compact manifold,  $n$ even,
the space $\Omega^{n/2}(M)$ of (complex) middle dimension forms
has a (complexif\/ied) inner product
\[
\langle\omega_1,\omega_2\rangle = \int_{M} \overline{\omega_1}
\wedge \star \omega_2.
\]
This inner product is unchanged under a conformal change of the
metric  and so,  its Hilbert space completion
$\mathcal{H}_0=L^2(M,\Omega^{n/2}(M))$ depends only on the
conformal class of the metric. $\mathcal{H}_0$~is by construction
a $C^\infty(M)$-module with $(f \omega)(p) = f(p)\omega(p)$ for
all $f \in C^\infty(M),$ $\omega\in \mathcal{H}_0,$ and $p \in M$.

If $M$ is a compact manifold without boundary then, the harmonic
forms (those in the kernel of $\Delta$) are precisely those in
$\kernel d \cap \kernel \delta.$   If $M$ is even dimensional with
dimension $n,$ the Hodge decomposition for middle-dimension forms
looks like
\[
\Omega^{n/2} (M) = \Delta(\Omega^{n/2}(M))\oplus H^{n/2} =
d(\Omega^{n/2-1}(M))\oplus \delta (\Omega^{n/2+1}(M)) \oplus
H^{n/2}.
\]
Here $H^{n/2}=\kernel_{n/2}\Delta.$ Thus $\mathcal{H}_0$ is the
direct sum of  $H^{n/2}$ and the images of $d$ and $\delta.$

The Hilbert space $\mathcal{H}$ is $\mathcal{H} = \mathcal{H}_0
\oplus H^{n/2},$ the direct sum of $\mathcal{H}_0$ with an extra
copy of the f\/inite dimensional Hilbert space of harmonic
middle-dimension forms on $M$.

For each $f\in \mathcal{A} =C^\infty(M)$ we consider the
multiplication operator on $\mathcal{H}_0,$ $f:\omega \mapsto
f\omega.$  These multiplication operators on $\mathcal{H}_0$ do
not preserve the subspace of harmonic forms.  Thus the extra copy
of $H^{n/2}$ in $\mathcal{H}$ is to preserve the notion of
$\mathbb{Z}_2$-graded Hilbert space.

The Hilbert space representation of $C^\infty(M)$ in $\mathcal{H}$
is given by $f \mapsto \pi(f)$ with $\pi(f)(\omega +h):= f\omega,$
for~all $\omega \in \mathcal{H}_0$ and $h \in H^{n/2}.$ Evidently,
$\pi(f)$ is a bounded operator on $\mathcal{H}.$ It is not
dif\/f\/icult to verify that $\pi$ is an involutive representation
of $C^\infty(M)$ in $\mathcal{H}.$ To simplify the notation we
write~$f$ instead of $\pi(f).$

Next we look at the $\mathbb{Z}/2$ grading.    Because the Hodge
star operator $\star$ acting on middle-dimension forms satisf\/ies
$\star^2 =  (-1)^{n/2},$ the operator
\[
\gamma_0 := (-1)^{\frac{(n/2)(n/2-1)}{2}} i^{n/2} \star
\]
is of square one giving a $\mathbb{Z}_2$-grading on
$\mathcal{H}_0.$ Since $\star^* = (-1)^{n/2}\star$ when acting on
middle forms we have
\begin{equation}
\label{gamma*} \gamma_0^* := (-1)^{\frac{(n/2)(n/2-1)}{2}}
(-i)^{n/2} \star^* = (-1)^{\frac{(n/2)(n/2+3)}{2}} i^{n/2} \star =
\gamma_0.
\end{equation}
It follows from the def\/inition of $\gamma_0$ that
\begin{lemma}
The operator $\gamma_0$ exchanges the subspaces
$d\big(\Omega^{n/2-1}\big)$ and $\delta\big(\Omega^{n/2+1}\big)$
and their closures.   As a consequence and because of
\eqref{gamma*}, $\gamma_0(H^{n/2}) = H^{n/2}.$
\end{lemma}

We def\/ine $\gamma:\mathcal{H}_0 \oplus H^{n/2} \to \mathcal{H}_0
\oplus H^{n/2}$ by $\gamma(\omega+h) := \gamma_0(\omega) -
\gamma_0 h.$ In this way, the extra copy of $H^{n/2}$ is endowed
with the opposite $\mathbb{Z}_2$-grading $-\gamma$:  $H^{n/2{\pm}}
=\{h \in H^{n/2} :\gamma h = \mp h\}$ and so $\mathcal{H}$ has the
$\mathbb{Z}_2$-grading given by $\mathcal{H}^+ = \mathcal{H}_0^+
\oplus {H^{n/2}}^+$ and $\mathcal{H}^- = \mathcal{H}_0^- \oplus
{H^{n/2}}^-.$

It is straightforward to verify that the operator $\gamma$
def\/ined on $\mathcal{H}$ satisf\/ies $\gamma=\gamma^*,$
$\gamma^2=1,$ and $\gamma f = f \gamma$ for all $f\in
C^\infty(M).$

The f\/irst step to def\/ine the operator $F$ is the following
observation
\begin{lemma}
For $\omega = d\beta + \delta \beta' \in
d(\Omega^{n/2-1}(M))\oplus \delta (\Omega^{n/2+1}(M))$,  the
operator $F_0\colon\mathcal{H}_0 \to \mathcal{H}_0 $ defined by
$F_0(d\beta + \delta \beta' ) :=  d\beta - \delta \beta'$ and
extended as zero over $H^{n/2}$ is a partial isometry such that
$F_0^2=1$ on $\mathcal{H}_0 \ominus H^{n/2}$, and $F_0$ is its own
formal adjoint operator.  Furthermore, $1 - F_0^2$ is the
orthogonal projection on the finite-dimensional Hilbert space of
middle-dimension harmonic forms.
\end{lemma}
The operator $F$ is def\/ined on $\mathcal{H} =
d(\Omega^{n/2-1}(M))\oplus \delta (\Omega^{n/2+1}(M)) \oplus
H^{n/2} \oplus H^{n/2}$ by{\samepage
\begin{equation}
\label{definitionF}
F = \begin{pmatrix} F_0 & 0 & 0 \\
                    0 & 0 & 1 \\ 0 & 1 & 0 \end{pmatrix}.
\end{equation}
From the previous lemma $F^* = F$ and $F^2 = 1$.}

For an even dimensional oriented compact manifold without
boundary, both the Hilbert space~$\mathcal{H}$ and the operator
$F$ are conformally invariant, thanks to the fact:  for a k-form
$\rho,$
\[
\widehat \delta \rho = e^{-(n-2(k-1))\eta}\,\delta
\,e^{(n-2k)\eta}\rho.
\]
It is not dif\/f\/icult to verify that $\gamma F = - F \gamma.$
Last, since each $[F,f]$ is a pseudodif\/ferential operator of
order~$-1$ for all $f \in \mathcal{A}=C^\infty(M),$ the operator
$[F,f]$ is a compact operator on $\mathcal{H}$ via Rellich's
theorem~\cite[p.~306]{Polaris}.

\subsection{Recovering the conformal structure}

Theorem~IV.4.2.c~of~\cite{Koran} states that the Fredholm module
$(\mathcal{H},F)$ uniquely determines the conformal structure of
$M$.   For that, Connes uses his trace theorem and the
noncommutative residue to recover the $L^n$-norm for exterior
$1$-forms over the manifold.

The f\/irst step is to consider instead of $d f$ its quantized
version $[F,f].$  Since $F$ is a pseudodif\/ferential operator of
order 0, $[F,f]$ is a pseudodif\/ferential operator of order $-1$
for all $f \in C^\infty(M^n),$ acting on the same vector bundle
$\Omega^{n/2} M$ as $F.$ The leading symbol of $F$ is given by
\[
\sigma_0(F)(x,\xi) =
|\xi|^{-2}\left(\varepsilon_{\tfrac{n}{2}-1}(\xi)
\iota_{\tfrac{n}{2}}(\xi)
                - \iota_{\tfrac{n}{2}+1}(\xi) \varepsilon_{\tfrac{n}{2}}(\xi)\right)
\]
for all $(x,\xi) \in T^* M$, $\xi \not= 0$.  Here
$\varepsilon_k(\xi)$ and $\iota_k(\xi)$ represent the exterior and
interior multiplication by the $1$-form $\xi$ on $k$-forms.  Note
how $\sigma_0(F)$ does not depend on $x\in M.$ The principal
symbol of $[F,f]$ is
\[
\sigma_{-1}([F,f])(x,\xi) = -i\sum_{k=1}^n \partial_{x^k}f
\,\partial_{\xi_k} (\sigma_0(F))
\]
which by the expression for $\sigma_0(F)(x,\xi)$ depends only on
the value on $x$ of the $1$-form $d f=\sum \partial_{x^k}f\,d
x^k.$  The details of these statements can be seen for example in
\cite{mio1}. Next for $f_i \in C^\infty(M)$ the operator
$(f_1[F,f_2])^n$ is a pseudodif\/ferential operator of order $-n.$

What Theorem~IV.4.2.c~\cite{Koran} shows is that
$\Wres(f_1[F,f_2])^n$ and $\int_M ||f_1 d f_2||^n \,d x,$ the
$L^n$-norm for $1$-forms, are proportional.

In the setting of spin Riemannian manifolds, for the algebra
$C^\infty(M)$ of smooth complex valued functions, the Hilbert
space is chosen to be $\mathcal{H} = L^2(M,\Sigma M),$ the square
integrable complex spinor f\/ields, and for $F$ one considers the
sign of the Dirac operator $D.$  Recently, B\"ar~\cite{Bar} showed
the following result.
\begin{theorem}
Let $M$ be a compact spin Riemannian manifold.  Let $g$ and $g'$
be Riemannian metrics on $M$ and let $(\mathcal{H}, \sign(D))$ and
$(\mathcal{H}',\sign(D'))$ be the corresponding Fredholm modules
of the algebra $C^\infty(M).$  Then $g$ and $g'$ are conformally
equivalent if and only if $(\mathcal{H}, \sign(D))$ and
$(\mathcal{H}',\sign(D'))$ are weakly unitarily equivalent.  That
is to say, there is a unitary isomorphism $U : \mathcal{H} \to
\mathcal{H}'$ such that $D'- UDU^{-1}$ is a compact operator and
for all $f \in C^\infty(M)$ and all $h \in \mathcal{H}$ one has
$U(f h) = f U(h).$
\end{theorem}

The idea is based on the commutativity of $U$ with the action of
$C^\infty(M)$ which implies that $U$ is induced by a (a.e.
invertible) section $\Psi$ of $L^\infty(M, \Hom(\Sigma M, \Sigma'
M)).$   The principal symbol of a~Dirac operator is given by
Clif\/ford multiplication with respect to the metric $g,$
$\sigma_D(\xi) = i c_g(\xi),\quad$  for all $\xi \in T^*M.$
Because of the relation $c_g(\xi)c_g(\eta) + c_g(\eta)c_g(\xi) =
-2g(\xi,\eta),$ for all $\xi,\eta \in T^*M,$ the principal symbol
of $\sign(D)$ is
\[
\sigma_{\sign(D)}(\xi) = \frac{i c_g(\xi)}{||\xi||_g},\qquad
\forall\, \xi \in T^*M \setminus \{0\}.
\]
Since $D'$ and $UDU^{-1}$ dif\/fer by a compact operator, they
have the same sign and thus
\[
\frac{c_{g'}(\xi)}{||\xi||_{g'}} = \Psi(x)
\frac{c_g(\xi)}{||\xi||_g}\Psi^{-1}(x)
\]
for all nonzero $\xi \in T^*M.$ Last
\begin{align*}
\frac{-2g'(\xi,\eta)}{||\xi||_{g'}||\eta||_{g'}} &=
\frac{c_{g'}(\xi)c_{g'}(\eta) + c_{g'}(\eta)c_{g'}(\xi)
}{||\xi||_{g'}||\eta||_{g'}}
\\
&=\Psi(x)\biggl(\frac{c_g(\xi)c_g(\eta)
   + c_{g}(\eta)c_{g}(\xi) }{||\xi||_{g}||\eta||_{g}}\biggr)\Psi^{-1}(x)
 =\frac{-2g(\xi,\eta)}{||\xi||_{g}||\eta||_{g}}
\end{align*}
since the term in the middle is a scalar.

It is important to recall here the result of Connes (see for
example \cite[p.~544]{Koran}) that says that one recovers the
metric distance between points in a connected manifold $(M,g)$
from the relation
\[
d(x,y) = \sup\{f(x) - f(y): f\in C^\infty(M) \hbox{~with~}
||[D,f]|| \leq 1\}.
\]
Note how with the stronger requirement $D' = UDU^{-1}$ (unitarily
equivalent) then $ ||[D,f]|| =  ||[D',f]||$ and thus $d=d'.$

If the conformal geometry of $(M,[g])$ is encoded in the Fredholm
module $(\mathcal{H},F,\gamma)$ over the algebra $C^\infty(M),$
then how can one extract the conformal geometry from this Fredholm
module?

One possibility is to use it to f\/ind conformal invariants
associated to a given conformal manifold, for example, as in  the
introduction.

\section{The noncommutative residue for manifolds with boundary}

\begin{remark}
Wodzicki: (see e.g.~\cite{Schrohe}) There is no non-zero trace on
the algebra of classical pseudodif\/ferential operators mod the
ideal of smoothing operators $\Psi^\infty(M)/\Psi^{-\infty}(M),$
whenever $M$ is noncompact or has a boundary.
\end{remark}

The noncommutative residue of Fedosov--Golse--Leichtnam--Schrohe
\cite{Fedosovetal} for manifolds with boundary is the unique (up
to a constant multiple) continuous trace for the operators in
Boutet de Monvel's algebra.  Roughly speaking, this noncommutative
residue acts on operators $A$ that are described by pairs of
symbols $\{a_i,a_b\}$ called interior and boundary symbol
respectively.  In case the manifold has empty boundary this
noncommutative residue coincides with the usual noncommutative
residue of Wodzicki, Guillemin, Adler, and Manin.

The setting for the noncommutative residue is given by a compact
manifold $M$ with boundary~$\partial M$ such that $M$ is embedded
in a compact manifold $\widetilde{M}$ without boundary, both~$M$
and~$\widetilde{M}$ of dimension $n >1.$  For $M$ we consider in a
boundary chart local coordinates given by $(x',x^n)$ with
$x'=(x^1,\dots,x^{n-1})$ coordinates for~$\partial M$ and $x^n$
the geodesic distance to $\partial M.$ It is important to mention
that the geodesic coordinate chosen for $x^n$ is only a technical
tool since the noncommutative residue is independent of the metric
and of local representations.

\subsection{Boutet de Monvel's sub-algebra of diagonal symbols}

In \cite{Boutet, Fedosovetal,Grubb}, and \cite{Schrohe} one can
f\/ind detailed introductions to Boutet de Monvel's calculus. The
operators in Boutet de Monvel's algebra we are interested in are
diagonal matrices of operators (endomorphisms) $A$ acting on
sections of vector bundles $E$ over $M$ and $E'$ over $\partial
M:$
\begin{align*}
A=
\begin{pmatrix} r^M P e^M + G & 0
\\
0&S
\end{pmatrix}:
\begin{matrix}
C^\infty(M,E)
\\
\oplus
\\
C^\infty(\partial M,E')
\end{matrix}
\to
\begin{matrix}
C^\infty(M,E)
\\
\oplus
\\
C^\infty(\partial M,E')
\end{matrix}.
\end{align*}
They are better described by a pair of symbols $(a_i, a_b)$ where
$a_i$ is called the interior symbol and $a_b$ is called the
boundary symbol. According to our needs, the characterization of
such an operator (or its symbol) of order $m$ is as follows.
\begin{itemize}\itemsep=0pt
\item[P.] The operator $P$ is a classical pseudodif\/ferential
operator of order $m$ on $\widetilde{M}.$  Furthermore, $P$~has
the so called transmission property.  This guarantees that the
composition of dif\/ferent elements remains inside the algebra.
Analytically, in local coordinates near $\partial M$ the
transmission property is given by
\[
\partial ^k_{x^n} \partial^\alpha_{\xi'} p_j (x',0,0,+1) = (-1)^{j-|\alpha|} \partial^k_{x^n} \partial^\alpha_{\xi'} p_j(x',0,0,-1),
\qquad \forall \,j,k,\alpha.
\]
Here $p_j$ is the homogeneous component of order $j$ in the symbol
expansion of the symbol~$p$ of~$P.$ Last, $e^M$ is the extension
by zero of functions (or sections) on $M$ to functions (or
sections) on $\widetilde{M}$ and $r^M$ is the restriction from
$\widetilde{M}$ to $M.$  The interior symbol $a_i$ of $A$ is
precisely $p.$
\end{itemize}
With $\mathcal{F}$ we denote the Fourier transform.  Also
\begin{align*}
H^+&=\{\mathcal{F}(\chi_{]0,\infty[}u) : u
\hbox{~is~a~rapidly~decreasing~function~on ~} \mathbb{R}\},
\\
H^-_0 &=\{\mathcal{F}((1-\chi_{]0,\infty[})u) : u
\hbox{~is~a~rapidly~decreasing~function~on ~} \mathbb{R}\},
\\
H^-&=H^-_0\oplus\{\hbox{all~polynomials}\}.
\end{align*}

The (diagonal) boundary symbol $a_b$ is given by a pair of symbols
$b$,~$s$ of operators~$G$,~$S$ parametrized by $T^* \partial
M\setminus \{0\}$ and the restriction of $p$ to the boundary.
\begin{itemize}
\item[G.] The operator $G$ is given by a singular green
operator-symbol $b(x',\xi',D_n)$ in the following way. For every
$l$ and f\/ixed $x'$, $\xi',$
\[
b_{l}(x',\xi',\xi_n,\eta_n) \in H^+ \hat \otimes_\pi H^-.
\]
With $\hat\otimes_\pi$ we denote Grothendieck's completion of the
algebraic tensor product. The ope\-rator $b(x',\xi',D_n) : H^+ \to
H^+$ is given by
\begin{gather*}
[b(x',\xi',D_n)h](\xi_n) = \Pi'_{\eta_n}
\bigl(b(x',\xi',\xi_n,\eta_n)h(\eta_n)\bigr) =\lim_{\eta_n \to
0^+} \mathcal{F}^{-1}(b(x',\xi',\xi_n,\cdot)h(\cdot))(\eta_n).
\end{gather*}

The operator $G$ described by this operator-symbol
$b(x',\xi',D_n)$ between functions on $[0,\infty[$ that are
rapidly decreasing at $\infty,$ def\/ines a trace class operator
on $L^2(R_+).$  The trace is given by
\[
\trace (G) (x',\xi') = \frac{1}{2\pi} \int b(x',\xi',\xi_n,\xi_n)
\,d\xi_n.
\]
Note that this is actually a symbol itself.

\item[S.] The operator $S$ is a classical pseudodif\/ferential
operator of order $m$ along the boundary.  It has values in
$\mathcal{L}(\mathbb{C}^k)$ and each component $s_j$ of its symbol
expansion $s$ acts by multiplication on $\mathbb{C}^k.$
\end{itemize}

The (diagonal) boundary symbol $a_b$ is then
\[
a_b(x',\xi',\xi_n,\eta_n)=
\begin{pmatrix}p(x',0,\xi',\xi_n) + b(x',\xi',\xi_n,\eta_n) & 0
\\
0 & s(x',\xi')
\end{pmatrix}
\]
with
\[
b(x',\xi',\xi_n,\eta_n) \sim \sum_{l=-\infty}^m
b_{l}(x',\xi',\xi_n,\eta_n) \qquad\hbox{and}\qquad s(x',\xi')\sim
\sum_{l=-\infty}^m s_l(x',\xi')
\]
where for $\lambda >0$
\[
b_{l}(x',\lambda \xi',\lambda \xi_n,\lambda \eta_n) = \lambda^{l}
b_{l}(x',\xi',\xi_n,\eta_n), \qquad s_{l}(x',\lambda \xi') =
\lambda^{l} s_{l}(x',\xi').
\]

By $\mathcal{B}_D^m(M)$ we denote the collection of all operators
of order $m$ with diagonal boundary symbol and by
$\mathcal{B}_D^\infty(M)$ the union of all the
$\mathcal{B}_D^m(M).$  The intersection over all orders $m$ of
$\mathcal{B}_D^m(M)$ is denoted $\mathcal{B}_D^{-\infty}(M).$ Last
$\mathcal{B}_D = \mathcal{B}_D^\infty(M) /
\mathcal{B}_D^{-\infty}(M).$

Given two operators $A_1$ and $A_2$ in $\mathcal{B}_D$ with
symbols $(a_{i1},a_{b1})$ and $(a_{i2},a_{b2}),$ with entries in
the boundary symbols $b_j$, $s_j,$ for $j=1,2,$ the composition is
again an operator in $\mathcal{B}_D$ with symbol $(a_{i},a_{b})$
where $a_i$ is the usual composition of symbols $a_i = a_{i1}
\circ a_{i2}.$  It also satisf\/ies the transmission property.

The resulting boundary symbol is of the form
\[
a_b = a_{b1} \circ' a_{b2} +
\begin{pmatrix}
L(p_{i1},p_{i2}) + p_{i1}^+ \circ' b_2 + b_1 \circ' p_{i2}^+ & 0
\\
0 &0
\end{pmatrix}.
\]
The symbol $\circ'$ denotes the usual composition of
pseudodif\/ferential symbols on the variables $(x',\xi').$  The
terms in the second summand represent the portion on the boundary
symbol coming from the interior symbols.  Here, we have hidden in
$a_{b1} \circ' a_{b2}$ the part corresponding to the restriction
to the boundary of the interior symbol.

The so called  ``left-over term'' $L(p_{i1},p_{i2}),$ ref\/lects
the particular way the pseudodif\/ferential operators $P_M = r^M P
e^M$ act on the manifold with boundary $M.$  If $P_1$ and $P_2$
are two pseudo\-dif\/ferential operators on $\widetilde M,$ the
dif\/ference $(P_1 P_2)_M-(P_1)_M (P_2)_M$ is a singular Green
operator with associated singular Green operator-symbol
$L(p_1,p_2).$ Since this left-over term need not be zero, we can
not reduce the diagonal sub-algebra by requesting $G=0$ in all the
operators.

As an example, and because they will be needed later on, let us
look at $L(f,q)$ and $L(p,f)$ where $p$ and $q$ are the symbols of
pseudodif\/ferential operators $P$ and $Q$ on $\widetilde{M},$
and $f \in C^\infty(\widetilde{M}),$ i.e.~$f$ represents the
pseudo\-dif\/ferential operator on $\widetilde{M}$ given
multiplication by $f.$ Among all the possible formulae for
$L(p,q)$ available in the literature we decided to use the one
provided in~\cite{Grubb}.

In Section~3 of \cite{Grubb} one can read an explicit expression
for $L(p,q)$ in which the ef\/fects of $p$ and $q$ are neatly
separated.  This expression uses singular Green operators $G^+(p)$
and $G^-(q)$ natural for the calculus in use (see Theorems~3.2
and~3.4 \cite{Grubb}).

We content ourselves by quoting a particular situation.  By
Theorem~3.4~\cite{Grubb}, $G^-(f)=0$ and by $(3.16)$~\cite{Grubb},
$L(p,f) = G^+(p)G^-(f),$ thus $L(p,f) = 0.$

Now, for $L(f,q)$ we must look at Theorem~3.5 \cite{Grubb}.  In
general,
\[
L(p,q) = G^+(p)G^-(q) +\sum_{0\leq m <
\hbox{\begin{tiny}~order~of~\end{tiny}} Q} K_m\gamma_m
\]
where the $K_m$ are operators obtained from symbols of a
particular type known as Poisson symbols.  By
$(3.35)$~\cite{Grubb}, $K_m =0$ when $p=f$ since it depends on
higher derivatives on $\xi_n.$  Since by Theorem~3.2~\cite{Grubb},
$G^+(f)=0$ we conclude that $L(f,p)=0$ as well.

\begin{lemma}
For every $f\in C^\infty(\widetilde M)$ and every
pseudodifferential operator $P$ on $\widetilde{M}$ with
symbol~$p$,  both left-over terms $L(f,p)$ and $L(p,f)$ vanish.
\end{lemma}
Last, the operator $p^+(x',\xi',D_n):H^+ \to H^+$ is induced from
the action (of the interior symbol) in the normal direction for
f\/ixed $(x',\xi').$  The only case we will be interested in are
those of the form $f^+ \circ' b_2$ where $f$ is a smooth function
on $\widetilde{M}.$  We will address them in \eqref{fplusb}.

\subsection{The noncommutative residue}

On $\mathbb{R}^n$ with coordinates $\xi_1,\dots,\xi_n$ we consider
the $(n-1)$-form
\[
\sigma = \sum_{j=1}^n (-1)^{j+1} \xi_j\, d\xi_1 \wedge \cdots
\wedge \widehat{d \xi_j}  \wedge \cdots \wedge d\xi_n,
\]
where the hat indicates this factor is omitted. Restricted to the
unit sphere $\mathbb{S}^{n-1},$  $\sigma$ gives the volume form on
$\mathbb{S}^{n-1}$ and in general $d\sigma = n \, d\xi_1 \wedge
\cdots \wedge d\xi_n.$  For a coordinate chart $U,$ the form $d
x_1 \wedge \cdots \wedge d x_n$ def\/ines an orientation on $U$
and induces the orientation $d\xi_1 \wedge \cdots \wedge d\xi_n$
on $\mathbb{R}^n.$

For a closed compact manifold $M$ without boundary, the
noncommutative residue is def\/ined as the unique trace (up to
constant multiples) on the algebra $\Psi^\infty/\Psi^{-\infty}$ of
classical pseudodif\/ferential operators mod the ideal of
smoothing operators.

The following is the main result of \cite{Fedosovetal}:
\begin{theorem}[Fedosov--Golse--Leichtnam--Schrohe]
Let $M$ be a manifold of dimen\-sion~$n$ with smooth boundary
$\partial M$, and let  $M\cup \partial M$ be embedded in a
connected manifold $\widetilde M$ of dimen\-sion~$n$. Let
\begin{align*}
A=
\begin{pmatrix} r^M P e^M + G & K
\\
T&S
\end{pmatrix}
\end{align*}
be an element in $\mathcal{B}^\infty(M)/\mathcal{B}^{-\infty}(M),$
with $\mathcal{B}^\infty(M)$ the algebra of all operators in
Boutet de Monvel's calculus (with integral order),
$\mathcal{B}^{-\infty}(M)$ the ideal of smoothing operators, and
let $p$,~$b$, and~$s$ denote the local symbols of $P$, $G$, and
$S$ respectively.  Then
\begin{align*}
\overline \Wres A =& \int_M \int_{S^{n-1}} \Trace_E p_{-n}(x,\xi)
\sigma(\xi) \,dx
\\&
+ 2\pi\int_{\partial M} \int_{S^{n-2}}
      \left\{\Trace_{E'} (\trace b_{-n})(x',\xi') + \Trace_{E'} s_{1-n}(x',\xi')\right\}\sigma'(\xi')\,dx',
\end{align*}
with $\sigma'$ the $n-2$ analog of $\sigma$, is the unique
continuous trace (up to constant multiples) on the algebra
$\mathcal{B}^\infty(M)/\mathcal{B}^{-\infty}(M)$.
\end{theorem}
This trace reduces to the noncommutative residue (of Adler, Manin,
Guillemin, and Wodzicki) in the case $\partial M = \varnothing,$
and it is independent of the Riemannian metric (eventually) chosen
on~$M.$

\section{On manifolds with boundary}
\label{sectionwithboundary}

In this section we present an extension of Theorem~1 in
\cite{Ugalde} to the setting of manifolds with boundary.  Let $M$
be a manifold with boundary $\partial M.$ Assume that the compact
manifold $M$ is embedded in a compact manifold $\widetilde M$
without boundary. Further we assume $\widetilde M$ to be oriented
which determines an orientation on $M$ and thus on $\partial M.$

For $P$ a pseudodif\/ferential operator acting on a vector bundle
$E$ over $\widetilde M$ with symbol $p$ having the transmission
property up to the boundary, $S$ a pseudodif\/ferential operator
acting on a vector bundle $E'$ over $\partial M$ with symbol $s,$
and for $f \in C^\infty(\widetilde M)$ we let $A(P,S)$ and $A(f)$
be the elements in Boutet de Monvel's algebra of diagonal elements
given by
\[
A(P,S) = \begin{pmatrix} r^M P e^M +0 & 0
\\
0 & S
\end{pmatrix},
\qquad A(f) =
\begin{pmatrix}
r^M f e^M +0 & 0
\\
0 & f|_{\partial M}
\end{pmatrix}.
\]
We study $\Wres(A(f_0)[A(P,S),A(f_1)][A(P,S),A(f_2)])$ for
functions $f_i \in C^\infty(\widetilde M).$

First of all, we must check that this product operator remains
inside the calculus in use.  It follows from
Proposition~2.7~of~\cite{Schrohe2}, which states that if two
operators satisfy the transmission property then their product
satisf\/ies the transmission property as well.

Since $L(f,p) =0 =L(p,f),$ for all $f \in C^\infty(\widetilde M)$
it follows that
\begin{align*}
&A(f_0)[A(P,S),A(f_1)][A(P,S),A(f_2)] =
\\
&
\begin{pmatrix}
r^M f_0[P,f_1][P,f_2] e^M + f_0^+ \circ'
L(\sigma([P,f_1]),\sigma([P,f_2]))& 0
\\
0 & f_0|_{\partial M} \circ' [S,f_1|_{\partial M}]' \circ'
[S,f_2|_{\partial M}]'
\end{pmatrix}
\end{align*}
where $\circ'$ represents the symbol composition with respect to
$(x',\xi').$   Here $\sigma([P,f_i])$ represents as usual the
symbol of the operator $[P,f_i].$  Using the def\/inition of
$\overline\Wres$ for manifolds with boundary we have
\begin{gather*}
\overline\Wres\bigl( A(f_0)[A(P,S),A(f_1)][A(P,S),A(f_2)] \bigr)
\\
\qquad{} = \int_M \int_{S^{n-1}} \Trace_E\left\{
\sigma_{-n}\bigl(f_0[P,f_1][P,f_2](x,\xi)\bigr) \right\}
\sigma(\xi) \,d x
\\
\qquad\quad{}+ 2\pi\int_{\partial M} \int_{S^{n-2}} \Trace_{E'}
\left\{ \sigma_{-(n-1)}\bigl( \trace\bigl\{ f_0^+ \circ'
L(\sigma([P,f_1]),\sigma([P,f_2])) \bigr\}(x',\xi') \bigr)
\right\}
\\
\qquad\qquad\qquad\qquad\quad {}+ \Trace_{E'} \left\{ \sigma_{-(n-1)}\left(
f_0|_{\partial M} \circ' [S,f_1|_{\partial M}]' \circ'
[S,f_2|_{\partial M}]' (x',\xi') \right) \right\} \sigma'(\xi')\,d
x'.
\end{gather*}

\subsection{A pair of bilinear functionals}

Mimicking the boundaryless case and following \cite{Wang} we
def\/ine:
\begin{definition}
\label{BandpartialB}
\[
B_{n,P}(f_1,f_2) := \int_{S^{n-1}}
\Trace_E\left\{\sigma_{-n}\bigl([P,f_1][P,f_2](x,\xi)\bigr)\right\}
\sigma(\xi),
\]
and
\begin{align*}
\partial B_{n,P,S}(f_1,f_2) := \int_{S^{n-2}}&
\Trace_{E'}
\left\{\sigma_{-(n-1)}\bigl(\trace\bigl\{L(\sigma(([P,f_1]),\sigma(([P,f_2]))\bigr\}
(x',\xi')\bigr)\right\}
\\
&{} + \Trace_{E'} \left\{\sigma_{-(n-1)}
\bigl(\bigl([S,f_1|_{\partial M}]' \circ'[S,f_2|_{\partial M}]\bigr)
(x',\xi') \bigr) \right\} \sigma'(\xi'),
\end{align*}
for all $f_i \in C^\infty(\widetilde M).$
\end{definition}
By def\/inition, both $B_{n,P}$ and $\partial B_{n,P,S}$ are
bilinear.  Since $f_0$ is independent of $\xi$ we have
\[
\int_{S^{n-1}} \Trace_E\left\{\sigma_{-n}\bigl(f_0[P,f_1][P,f_2]
(x,\xi)\bigr)\right\} \sigma(\xi)  = f_0 B_{n,P}(f_1,f_2).
\]
The computations done in \cite{MioCMP} with the symbol expansions
for the case of empty boundary are also valid here.  In particular
we have in given local coordinates the explicit expression
\[
B_{n,P}(f_1,f_2) = \sum \frac{D^\beta_x(f_1)
D^{\alpha''+\delta}_x(f_2)}{\alpha'!\alpha''!\beta!\delta!}
\int_{S^{n-1}} \Trace\biggl\{
\partial^{\alpha'+\alpha''+\beta}_\xi(\sigma^P_{k-i}) \partial^{\delta}_\xi(D^{\alpha'}_x(\sigma^P_{k-j})) \biggr\}\sigma(\xi)
\]
with the sum taken over $|\alpha'|+|\alpha''|+|\beta|+|\delta|+i+j
= n+2k,  |\beta|\geq 1$, and $|\delta|\geq 1.$  It shows that
$B_{n,P}(f_1,f_2)$ is dif\/ferential in $f_1$ and $f_2.$
Evidently it is possible to obtain a similar expression for the
summand in $\partial B_{n,P,S}$ corresponding to $S$ replacing $n$
by $n-1$ and $x$ by $x'.$

In p.~25 of \cite{Fedosovetal} we can read an expression for the
degree $-(n-1)$ component of the operator-symbol $\trace(c)$ with
$c=p^+\circ' b.$  It is given by
\begin{gather}
\label{fplusb} \sigma_{-(n-1)} \bigl(\trace c(x',\xi')\bigr)
\!\sim\!\! \sum_{j=0}^\infty \frac{i^j}{j!} \Pi'_{\xi_n}
\bigl\{\sigma_{-n}\bigl(
\partial^j_{\xi_n}[\partial^j_{x_n} p(x',0,\xi',\xi_n)
\circ'
b(x',\xi',\xi_n,\eta_n)]\bigr)|_{\eta_n=\xi_n}\bigr\}.\!\!\!\!
\end{gather}
Thus, since $f_0$ is independent of $\xi,$
\begin{gather*}
\trace\bigl\{\sigma_{-(n-1)}\bigl(f_0^+ \circ'
L(\sigma([P,f_1]),\sigma([P,f_2]))\bigr)\bigr\}
\\
\qquad{}=
f(x',0)\trace\{\sigma_{-(n-1)}(L(\sigma([P,f_1]),\sigma([P,f_2]))(x',\xi')\}.
\end{gather*}
and it follows that
\begin{gather*}
\int_{S^{n-2}} \Trace
\left\{\sigma_{-(n-1)}\bigl(\trace\bigl\{f_0^+ \circ'
L(\sigma([P,f_1]),\sigma([P,f_2]))\bigr\}
(x',\xi')\bigr)\right\}\sigma'(\xi')
\\
\qquad{}=f(x',0) \int_{S^{n-2}} \Trace
\left\{\sigma_{-(n-1)}\bigl(\trace\bigl\{L(\sigma([P,f_1]),\sigma([P,f_2]))\bigr\}
(x',\xi')\bigr)\right\}\sigma'(\xi'),
\end{gather*}
for all $f_i \in C^\infty(\widetilde M).$  In this way
\begin{gather}
\nonumber \overline \Wres\bigl(A(f_0)[A(P,S) , A(f_1)] [A(P,S) ,
A(f_2)]\bigr)
\\
\qquad{}= \int_M f_0 B_{n,P}(f_1,f_2\bigr) \, d x +  2\pi
\int_{\partial M} f_0|_{\partial M}\partial
B_{n,P,S}(f_1|_{\partial M},f_2|_{\partial M}) \,
dx'.\label{foruniqueness}
\end{gather}
\begin{lemma}
The functionals $B_{n,P}$ and $\partial B_{n,P,S}$ are bilinear
and symmetric.
\end{lemma}
\begin{proof}
The symmetry  of both $B_{n,P}$ and $\partial B_{n,P,S}$ is not
evident from the expressions above.   For~$B_{n,P}$ it was
obtained  in \cite{mio1} in the boundaryless case from the trace
property of $\Wres.$  Because it shares the same local expression
both for empty and non-empty boundary we have that $B_{n,P}$ is
symmetric.

For $\partial B_{n,P,S}$ we are going to exploit the linearity and
the trace property of the noncommutative residue.  Denote
$\overline f =A(f)$ and $\overline P = A(P,S).$ Using that
$\overline{f_1} \, \overline{f_2} = \overline{f_2} \,
\overline{f_1}$ for all $f_i \in C^\infty(\widetilde M)$ and  the
trace property of the noncommutative residue we have that all of
$\overline\Wres\bigl(\overline{f_0} \, \overline{f_2} \, \overline
P \, \overline P \, \overline{f_1}
                     - \overline{f_0} \, \overline{f_1} \, \overline P \, \overline P \, \overline{f_2}\bigr),$
$\overline\Wres\bigl(\overline{f_0} \, \overline{f_2} \, \overline
P \, \overline P \, \overline{f_1}
                     - \overline{f_0} \, \overline{f_1} \, \overline P \, \overline P \, \overline{f_2}\bigr),$ and
$\overline\Wres\bigl(\overline{f_0} \, \overline P \,
\overline{f_2} \, \overline P \, \overline{f_1}
                     - \overline{f_1} \, \overline{f_0} \, \overline P \, \overline{f_2} \, \overline P\bigr)$ vanish.
In this way
\begin{gather*}
\overline\Wres(\overline{f_0}[\overline
P,\overline{f_1}][\overline P,\overline{f_2}] -
\overline{f_0}[\overline P,\overline{f_2}][\overline
P,\overline{f_1}])
\\
\qquad{}=\overline\Wres(\overline{f_0} \, \overline P \,
\overline{f_1} \, \overline P \, \overline{f_2} - \overline{f_0}
\, \overline P \, \overline{f_1} \, \overline{f_2} \, \overline P
- \overline{f_0} \, \overline{f_1} \, \overline P \, \overline P
\, \overline{f_2} + \overline{f_0} \, \overline{f_1} \, \overline
P \, \overline{f_2} \, \overline P
\\
\quad\qquad {}- \overline{f_0} \, \overline P \, \overline{f_2} \,
\overline P \, \overline{f_1} + \overline{f_0} \,\overline P \,
\overline{f_2} \, \overline{f_1} \, \overline P + \overline{f_0}
\, \overline{f_2} \, \overline P \, \overline P \, \overline{f_1}
- \overline{f_0} \, \overline{f_2} \, \overline P \,
\overline{f_1} \, \overline P)
\\
\qquad{}=\overline\Wres(\overline{f_0} \, \overline P \,
\overline{f_1} \, \overline P \, \overline f_2 + \overline{f_0} \,
\overline{f_1} \, \overline P \, \overline{f_2} \, \overline P -
\overline{f_0} \, \overline P \, \overline{f_2} \, \overline P \,
\overline{f_1} - \overline{f_0} \, \overline{f_2} \, \overline P
\, \overline{f_1} \, \overline P)=0.
\end{gather*}
Hence
\begin{gather*}
\int_M f_0 B_{n,P}(f_1,f_2) + 2\pi \int_{\partial M}
f_0|_{\partial M} \partial B_{n,P,S}(f_1,f_2)
\\
\qquad= \int_M f_0 B_{n,P}(f_2,f_1) + 2\pi \int_{\partial M}
f_0|_{\partial M} \partial B_{n,P,S}(f_2,f_1), \qquad \forall\,
f_i \in C^\infty(\widetilde M).
\end{gather*}
Since $B_{n,P}(f_1,f_2)$ is symmetric
\[
\int_{\partial M} f_0|_{\partial M} \overline B_S(f_1,f_2) =
\int_{\partial M} f_0|_{\partial M} \overline B_S(f_2,f_1),\quad
\forall f_0 \in C^\infty(\widetilde M)
\]
and the result follows from the arbitrariness of $f_0.$
\end{proof}
\begin{lemma}
$\partial B_{n,P,S}(f_1,f_2)$ is differential on $f_1$ and $f_2.$
\end{lemma}
\begin{proof}
We denote, to simplify the notation, $P_1=[P,f_1]$ and
$P_2=[P,f_2]$ with symbols~$p_1$ and~$p_2$ respectively.  In p.~27
of \cite{Fedosovetal} we can read the following
\begin{gather*}
\trace\bigl\{L(p_1,p_2)\bigr\} (x',\xi')
\\
\qquad{}= \sum_{j,k=0}^\infty \frac{(-i)^{j+k+1}}{(j+k+1)!}
\Pi'_{\xi_n}\bigl(\partial_{x_n}^j \partial^k_{\xi_n}
\Pi^+_{\xi_n}(p_1)(x',0,\xi',\xi_n) \circ'
                         \partial_{x_n}^{j+1} \partial^k_{\xi_n} \Pi_{\xi_n}^+(p_2)(x',0,\xi',\xi_n) \bigr),
\end{gather*}
with $\Pi^+_{\xi_n}(s)(\cdot)$ the projection of the symbol $s$ on
$H^+.$  The subscript in $\Pi^+$ indicates the variable it is
acting on.   From \cite{mio1} we know
\[
\sigma_{-k}([P,f]) = \sum_{|\beta|=1}^k
\frac{1}{\beta!}D^\beta_x(f)
\partial^\beta_\xi(\sigma^P_{-(k-|b|)})
\]
thus
\begin{align*}
\Pi_{\xi_n}^+(\sigma_{-k}([P,f]))(x',0,\xi',\xi_n) &=
\sum_{|\beta|=1}^k \frac{1}{\beta!}\Pi_{\xi_n}^+\bigl(D^\beta_x(f)
\partial^\beta_\xi(\sigma^P_{-(k-|b|)})\bigr)(x',0,\xi',\xi_n)
\\
&= \sum_{|\beta|=1}^k
\frac{1}{\beta!}D^\beta_x(f)(x',0)\Pi_{\xi_n}^+\bigl(\partial^\beta_\xi(\sigma^P_{-(k-|b|)})\bigr)(x',0,\xi',\xi_n).
\end{align*}
Since any $\partial^j_{x_n} f$ factors out of $\Pi'_{\xi_n}$ we
conclude the result.
\end{proof}
%

\subsection[Conformal invariance of $B_{n,P}$ and $\partial B_{n,P,S}$]{Conformal invariance of $\boldsymbol{B_{n,P}}$ and $\boldsymbol{\partial B_{n,P,S}}$}

If we further assume  Riemannian structures $(M,g)$ and
$(\widetilde M, \widetilde g)$ such that $g$ coincides with
$\widetilde g$ restricted to $M$ then, a conformal rescaling of
$g$ corresponds to a conformal rescaling of $\widetilde{g}$ (by an
appropriate extension of the conformal factor) and a conformal
rescaling of $\widetilde{g}$ can be restricted to a conformal
rescaling of~$g.$  We obtain
\begin{lemma}
\label{conformalrescaling} Assume that $P$ and $S$ are such that
$[P,f_1][P,f_2]$ and $[S,f_1|_{\partial M}][S,f_2|_{\partial M}]$
are conformally invariant for all $f_i \in C^\infty(\widetilde
M).$  Then
\[
\widehat{B_{n,P}}(f_1,f_2)(x)= e^{-2n\eta(x)}B_{n,P}(f_1,f_2)(x)
\]
and
\[
\widehat{\partial B_{n,P,S}}(f_1,f_2)(x') =
e^{-2(n-1)\eta(x',0)}\partial B_{n,P,S}(f_1,f_2)(x').
\]
\end{lemma}
\begin{proof}
We want to exploit the independence of $\overline \Wres$ of local
representations.  We have
\begin{gather*}
\overline\Wres(A(f_0)[A(P,S),A(f_1)][A(P,S),A(f_2)])
\\
\qquad{}= \int_M f_0 B_{n,P}(f_1,f_2)\,d x + 2\pi \int_{\partial
M} f_0|_{\partial M} \partial B_{n,P,S}(f_1|_{\partial
M},f_2|_{\partial M})\, d x'
\\
\qquad{}= \int_M f_0 \widehat{B_{n,P}}(f_1,f_2)\,\widehat{d x} +
2\pi \int_{\partial M} f_0|_{\partial M} \widehat{\partial
B_{n,P,S}}(f_1|_{\partial M},f_2|_{\partial M})\, \widehat{d x'}
\\
\qquad{}= \int_M f_0 e^{-2n\eta}\widehat{B_{n,P}}(f_1,f_2)\, d x +
2\pi \int_{\partial M} f_0|_{\partial M} e^{-2(n-1)\eta|_{\partial
M}}\widehat{\partial B_{n,P,S}}(f_1|_{\partial M},f_2|_{\partial
M})\, d x',
\end{gather*}
where we use $\,\widehat{\,}\,$ to represent quantities computed
with respect to the conformal metric $\widehat g = e^{2\eta}g.$
In particular
\[
\int_M f_0(x) B_{n,P}(f_1,f_2)(x)\, d x = \int_M f_0(x)
e^{-2n\eta(x)}\widehat{B_{n,P}}(f_1,f_2)(x)\, d x
\]
for all $f_i \in C^\infty(\widetilde M)$ with $f_0|_{\partial
M}=0.$  Thus $B_{n,P}(f_1,f_2)(x) =
e^{-2n\eta(x)}\widehat{B_{n,P}}(f_1,f_2)(x)$ for all $x \in M.$
It follows that
\[
\int_{\partial M} f_0(x',0)
\widehat{B_{n,P}}(f_1,f_2)(x')\,\widehat{d x'} = \int_{\partial M}
f_0(x',0) e^{-2n\eta(x',0)}\widehat{B_{n,P}}(f_1,f_2)(x')\, d x'
\]
for all $f_i \in C^\infty(\widetilde M).$ The result follows from
the arbitrariness of $f_0.$
\end{proof}
\begin{remark}
Note how the same reasoning in the proof above can be used to show
the uniqueness of $B_{n,P}$ and $\partial B_{n,P,S}$ satisfying
\eqref{foruniqueness}.
\end{remark}

Summarizing this section we have
\begin{theorem}
\label{Theorem1boundary} Let $M$ be a compact manifold of
dimension $n$ and with boundary $\partial M.$  Assume that $M$ is
embedded in a compact oriented manifold $\widetilde M$ without
boundary.  Further assume Riemannian structures
$(\widetilde{M},\widetilde{g})$ and $(M,g)$ such that $g$
coincides with $\widetilde{g}$ restricted to $M.$  Let $P$ be
a~pseudodifferential operator acting on a vector bundle $E$ over
$\widetilde M$ having the transmission property up to $\partial
M,$ let $S$ be a pseudodifferential operator acting on a vector
bundle $E'$ over $\partial M,$ such that $[P,f_1][P,f_2]$ and
$[S,f_1|_{\partial M}][S,f_2|_{\partial M}]$ are conformally
invariant for all $f_i \in C^\infty(\widetilde M).$  Then
$B_{n,P}$ and $\partial B_{n,P,S}$ given in Definition~{\rm
\ref{BandpartialB}} are conformally invariant in the sense of
Lemma~{\rm \ref{conformalrescaling}}. Furthermore, both $B_{n,P}$
and $\partial B_{n,P,S}$ are symmetric, bilinear differential
functionals uniquely determined by the
relation~\eqref{foruniqueness}.
\end{theorem}
%

\section{On even-dimensional manifolds with boundary}

Up to this point, we have a generalization of Theorem~1 in
\cite{mio1} to manifolds with boundary in the setting described
above.  Next, we want to state a generalization of Theorem~2 in
\cite{mio1} to this context.  In order to do it, we consider the
Fredholm module $(\mathcal{H},F)$ now associated to the
even-dimensional manifold without boundary $\widetilde{M}.$

\subsection[The symbol of $F$ and the transmission property]{The symbol of $\boldsymbol{F}$ and the transmission property}

If $\omega = d \beta + \delta \beta' \in
d(\Omega^{n/2-1}(\widetilde{M}))\oplus \delta
(\Omega^{n/2+1}(\widetilde{M}))$ then
\begin{gather*}
\Delta F_0 (d\beta + \delta \beta') = \Delta (d\beta
-\delta\beta') = d\delta d \beta - \delta d \delta\beta' =
F_0(d(\delta d \beta) + \delta (d\delta \beta')) = F_0 \Delta
(d\beta + \delta\beta').
\end{gather*}
It follows
\begin{lemma}
For an oriented compact manifold without boundary $\widetilde{M}$
and of even dimension $n,$ the relation $F_0 \Delta = \Delta F_0 =
d\delta - \delta d$ holds on
$d(\Omega^{n/2-1}(\widetilde{M}))\oplus \delta
(\Omega^{n/2+1}(\widetilde{M})).$
\end{lemma}

To be able to use a given pseudodif\/ferential operator in the
machinery of the noncommutative residue for manifolds with
boundary, it is essential for the operator to enjoy the
transmission property up to the boundary of $M.$

Because we are interested in $F$ acting on the orthogonal
complement of the harmonic forms on $\widetilde M,$ we abuse of
the notation  and use freely $F$ for $F_0.$  From the relation
$\Delta F = d\delta - \delta d$ and the formula for the total
symbol of the product of pseudodif\/ferential operators we can
compute the symbol expansion of $F.$  First we note that $F$ is a
pseudodif\/ferential operator of order $0.$

We know $\sigma(\Delta F) = \sigma(d\delta - \delta d),$ thus the
formula for the total symbol of the product of two
pseudodif\/ferential operators implies
\begin{align*}
\sigma_2^{d\delta-\delta d} + \sigma_1^{d\delta-\delta d} +
\sigma_0^{d\delta-\delta d} &= \sigma(d\delta - \delta d) =
\sigma(\Delta F) \sim \sum\frac{1}{\alpha!}\partial^\alpha_\xi
\sigma(\Delta) D^\alpha_x (\sigma(F))
\\
&\sim \sum \frac{1}{\alpha !} \partial^\alpha_\xi(\sigma_2^\Delta
+ \sigma_1^\Delta
  + \sigma_0^\Delta) D^\alpha_x(\sigma_0^F + \sigma_{-1}^F + \sigma_{-2}^F + \cdots).
\end{align*}
Expanding the right hand side into sum of terms with the same
homogeneity we conclude:

\begin{lemma}
\label{lemmasymbolofS} In any given system of local charts, we can
express the total symbol of $F$, $\sigma(F) \sim \sigma_0^F +
\sigma_{-1}^F + \cdots$ in a recursive way by the formulae:
\begin{gather*}
\sigma_0^F = (\sigma_2^\Delta)^{-1} \sigma_2^{d\delta-\delta d},
\qquad \sigma_{-1}^F =
(\sigma_2^\Delta)^{-1}\Bigl(\sigma_1^{d\delta-\delta d} -
\sigma_1^\Delta \sigma_0^F
 - \sum_{|\alpha|=1} \partial^\alpha_\xi(\sigma_2^\Delta) D^\alpha_x(\sigma_0^F) \Bigr),
\\
\sigma_{-2}^F =
(\sigma_2^\Delta)^{-1}\Bigl(\sigma_0^{d\delta-\delta d} -
\sigma_1^\Delta \sigma_{-1}^F - \sigma_0^\Delta \sigma_0^F
\\
\phantom{\sigma_{-2}^F =}{}
   - \sum_{|\alpha|=1}\left(\partial^\alpha_\xi(\sigma_2^\Delta) D^\alpha_x(\sigma_{-1}^F)
   + \partial^\alpha_\xi(\sigma_1^\Delta) D^\alpha_x(\sigma_0^F)\right)
   - \sum_{|\alpha|=2} \frac{1}{\alpha!} \partial^\alpha_\xi(\sigma_2^\Delta)D^\alpha_x(\sigma_0^F) \Bigr),
\\
\sigma_{-r}^F = -(\sigma_2^\Delta)^{-1}\Bigl( \sigma_1^\Delta
\sigma_{-r+1}^F + \sigma_0^\Delta \sigma_{-r+2}^F
   +\sum_{|\alpha|=1}\partial^\alpha_\xi(\sigma_2^\Delta) D^\alpha_x(\sigma_{-r+1}^F)
\\
\phantom{\sigma_{-r}^F =}{}
   +\sum_{|\alpha|=1}\partial^\alpha_\xi(\sigma_1^\Delta) D^\alpha_x(\sigma_{-r+2}^F)
   + \sum_{|\alpha|=2} \frac{1}{\alpha!} \partial^\alpha_\xi(\sigma_2^\Delta)D^\alpha_x(\sigma_{-r+2}^F)\Bigr),
\end{gather*}
for every $r \geq 3.$
\end{lemma}

Lemma~2.4~of~\cite{Schrohe2} states that all symbols which are
polynomial in $\xi$ have the transmission property.  Thus both
$\Delta$ and $d\delta - \delta d$ have the transmission property.
Proposition~2.7 in the same reference states that if two operators
satisfy the transmission property then their products, all their
derivatives, and their parametrizes satisfy the transmission
property as well.  Furthermore, the same result also states that
it is enough to check that each homogeneous component of the
symbol expansion has the transmission property to conclude that
the full symbol has the transmission property.

By Lemma~\ref{lemmasymbolofS}, each homogeneous component
$\sigma_{-k}^F$ in the symbol expansion of $F$ is given in terms
of derivatives of the homogeneous components of $\Delta,$
$d\delta-\delta d,$ $\sigma_0^F, \dots, \sigma_{-k+1}^F,$ and
$\sigma_2(\Delta)^{-1}.$  By Lemma~2.4 and Proposition~2.7
of~\cite{Schrohe2} it follows that
\begin{lemma}
The operator $F$ satisfies the transmission property.
\end{lemma}
%

\subsection[${\rm Res} ...$]{$\boldsymbol{\overline\Wres(A(f_0)[A(F,0),A(f_1)][A(F,0),A(f_2)])}$}

For $F$ given in \eqref{definitionF} now for the manifold
$\widetilde{M},$ and for $f \in C^\infty(\widetilde M)$ we let
$A(F,0)$ and $A(f)$ be the elements in Boutet de Monvel's algebra
of diagonal elements given by
\[
\overline F=A(F,0) = \begin{pmatrix} r^M F e^M +0 & 0
\\
0 & 0
\end{pmatrix},
\qquad \overline f=A(f) =
\begin{pmatrix}
r^M f e^M +0 & 0
\\
0 & f|_{\partial M}
\end{pmatrix}.
\]

Since $L(f,\sigma(F)) =0 =L(\sigma(F),f),$ it follows that
\[
\overline{f_0} [\overline F , \overline{f_1}] [\overline F ,
\overline{f_2}] =
\begin{pmatrix}
r^M f_0[F,f_1][F,f_2] e^M + f_0^+ \circ'
L(\sigma([F,f_1]),\sigma([F,f_2])) & 0
\\
0 & 0
\end{pmatrix},
\]
where $\circ'$ represents the symbol composition with respect to
$(x',\xi').$ As before we def\/ine:
\begin{definition}
\[
B_n(f_1,f_2) := \int_{S^{n-1}}
\Trace\left\{\sigma_{-n}\bigl([F,f_1][F,f_2] (x,\xi)\bigr)\right\}
\sigma(\xi),
\]
and
\[
\partial B_n(f_1,f_2)
= \int_{S^{n-2}} \Trace
\sigma_{-(n-1)}\left\{\bigl(\trace\bigl\{L(\sigma([F,f_1]),\sigma([F,f_2]))
(x',\xi')\bigr\}\bigr)\right\}\sigma'(\xi'),
\]
for all $f_i \in C^\infty(\widetilde M).$
\end{definition}
As in Section~\ref{sectionwithboundary}
\begin{theorem}
\label{Theorem2boundary} Both differential functionals $B_n$ and
$\partial B_n$ are bilinear, symmetric, conformal invariant in the
sense
\[
\widehat{B_n}(f_1,f_2)(x)= e^{-2n\eta(x)}B_n(f_1,f_2)(x)
\]
and
\[
\widehat{\partial B_n(f_1,f_2)(x')} =
e^{-2(n-1)\eta(x',0)}\partial B_n(f_1,f_2)(x')
\]
for a conformal change of the metric $\widehat g = e^{2\eta}g,$
and are uniquely determined by the relation
\begin{gather*}
\overline \Wres\bigl(A(f_0)[A(F,0) , A(f_1)] [A(F,0) ,
A(f_2)]\bigr)
\\
\qquad= \int_M f_0 B_{n,P}(f_1,f_2) \, d x +  2\pi \int_{\partial
M} f_0|_{\partial M}\partial B_{n,P,S}(f_1,f_2) \, dx'.
\end{gather*}
\end{theorem}
\begin{remark}
Even though both bilinear functionals $B_n$ and $\partial B_n$ are
acting on $C^\infty(\widetilde M),$ they depend on the particular
embedding of the compact manifold $M$ into $\widetilde M,$ and
thus, they can be def\/ined on $C^\infty(M)$ by considering an
extension of $f\in C^\infty(M)$ to $C^\infty(\widetilde M).$
\end{remark}
\begin{remark}
In case $M$ is odd dimensional, all results from the f\/irst part
of these notes are valid on the compact even dimensional manifold
without boundary $\partial M.$  In this way, we can consider the
commutative algebra $\mathcal{A}=C^\infty(\partial M)$  and the
Fredholm module associated to the manifold~$\partial M.$ For $F$
given in \eqref{definitionF} on the manifold~$\partial M,$ we
could look at
\[
A(P,F) = \begin{pmatrix} P & 0
\\
0 & F
\end{pmatrix}
\]
and try to study $\Wres(A(f_0)[A(P,F),A(f_1)][A(P,F),A(f_2)])$ for
functions $f_i \in C^\infty(\widetilde M).$  The trivial choice
$P=0$ will produce $B_{n,P} =0$ and $\partial B_{n,0,F} =
B_{n-1}.$  It is an open problem to search for a companion $P$ for
$F$ that will produce more interesting results in the odd
dimensional case.
\end{remark}
%

\subsection*{Acknowledgements}

This research is supported by Vicerrector\'{i}a de Investigaci\'on
de la Universidad de Costa Rica and Centro de Investigaciones
Matem\'aticas y Meta-matem\'aticas.  The material extends a talk
presented in May 2007 at the Midwest Geometry Conference held at
the University of Iowa in honor of Thomas P.~Branson.

The referees' suggestions improved to a great extent the
presentation of this material.  One of the referees pointed the
author towards \cite{Grubb} which provided a clearer understanding
of Boutet de Monvel's calculus.  In particular, the formulae used
for $L(p,q)$ resulted in a signif\/icant simplif\/ication of the
treatment of the subject.

\pdfbookmark[1]{References}{ref}
\LastPageEnding

\end{document}